\documentclass{article}

\usepackage{arxiv}

\usepackage[utf8]{inputenc} 
\usepackage[T1]{fontenc}    
\usepackage{hyperref}       
\usepackage{url}            
\usepackage{booktabs}       
\usepackage{amsfonts}       
\usepackage{nicefrac}       
\usepackage{microtype}      
\usepackage{lipsum}		
\usepackage{graphicx}
\usepackage{natbib}
\usepackage{doi}

\usepackage{amsmath,amsthm,amsfonts,amssymb,amsopn,cancel}
\usepackage{graphicx}
\usepackage{hyperref}
\usepackage{latexsym}
\usepackage{caption}
\usepackage{subcaption}

\title{A Convolutional Dispersion Relation Preserving Scheme for the Acoustic Wave Equation}
\author{{Oded Ovadia}\thanks{Corresponding author} \\
	Department of Applied Mathematics\\
	Tel Aviv University\\
	Tel Aviv 69978, Israel \\
	\texttt{odedovadia@mail.tau.ac.il}
	\And
	{Adar Kahana} \\
	Department of Applied Mathematics\\
	Tel Aviv University\\
	Tel Aviv 69978, Israel 
	\And
	{Eli Turkel} \\
	Department of Applied Mathematics\\
	Tel Aviv University\\
	Tel Aviv 69978, Israel 
	}

\begin{document}
\maketitle

	\begin{abstract}
	    We propose an accurate numerical scheme for approximating the solution of the two dimensional acoustic wave problem.
		We use machine learning to find a stencil suitable even in the presence of high wavenumbers.
		The proposed scheme incorporates physically informed elements from the field of optimized numerical schemes into a convolutional optimization machine learning algorithm.
	\end{abstract}
	
	\keywords{Numerical methods \and Optimization \and Dispersion \and Machine learning \and Physics-informed}

    \section{Introduction}

We propose a method to solve the wave propagation problem in a homogeneous domain. We aim to find a stencil that produces accurate results even in the setting of low resolution (coarse mesh) and high frequencies. We focus on the two-dimensional problem in the time domain (two spatial and one time coordinate). This problem has applications in various fields such as fluid dynamics, acoustics, electromagnetic problems, and more \cite{CFD}. The proposed method is based on a Physically-Informed (PI) Machine-Learning (ML) approach \cite{PINNs}.

Finite-difference (FD) schemes are a common numerical method for solving hyperbolic partial differential equations (PDEs) \cite{hyper1, hyper2, hyper3}. The physical domain is discretized via a mesh. The solution of the PDE is approximated using predefined numerical formulae on this mesh, such that each point is calculated using adjacent points in the mesh. This set of formulae defines an FD stencil. The main differences between various FD schemes is the choice of the stencil. Some stencils require more adjacent points for the calculation and others use more complex formulae. The different schemes aim to achieve accurate results compared to the analytic solution.

The field of ML is accelerating rapidly in the applied mathematics community. Many PDE related problems can be formulated as data-driven problems. ML models aim to map input data with corresponding outputs. They excel in non-linear problems and are usually able to fit the problem with high accuracy and quick convergence. However, one needs to either choose or design a specific ML method to fit the problem, otherwise it might not converge. By combining the classical methods with advanced ideas from the field of ML (such as PI neural networks) we achieve better results, as presented in this work.

The dispersion relation is a property of the wave problem that is based on three physical attributes: the wavenumber (spatial oscillations), frequency (temporal oscillations), and propagation velocity. When numerically solving a PDE using an FD scheme, the physical dispersion relation is replaced by another induced by the numerical scheme. The numerical dispersion relation could be quite different from the physical one. As a result, the FD approximation deviates from the analytic solution of the physical problem. This phenomenon is called numerical dispersion.

Numerical dispersion is dependent on the physical parameters of the problem, the choice of discretization, and the choice of FD method. Generally, numerical dispersion is more significant when dealing with high wavenumbers and negligible with low wavenumbers. Different numerical schemes induce different numerical dispersion relations based on their stencils. There are four common approaches to deal with numerical dispersion:

\begin{itemize}
	\item \textbf{Physical approach:} if possible, one can limit the wavenumber and solve only for low wavenumbers. In this case, the numerical dispersion has less impact on the accuracy, but valuable data is discarded.
	
	\item \textbf{Refining the mesh:} by solving on more nodes in the mesh, the numerical dispersion can be reduced. However, increasing the size of the mesh could be very expensive in terms of computational cost. 
	
	\item \textbf{High accuracy:} this approach uses a Taylor series to closely approximate the low wavenumbers. By increasing the order of accuracy of the method numerical dispersion could be mitigated \cite{high_order}. 
	
	\item \textbf{Using dispersion relation preserving (DRP) schemes:} this class of numerical schemes aims to create a stencil that produces accurate results on low resolution mesh with high wavenumbers. The method presented in the work is a DRP scheme.
\end{itemize}

The main idea of DRP schemes is to replace classical FD stencils with different ones for which the numerical dispersion relation represents the physical one more faithfully in the setting of high wavenubmers or using a coarse mesh. There is a variety of ways to obtain such stencils, mainly using different optimization techniques. The first DRP scheme was proposed in the seminal paper by Tam and Webb \cite{tam_and_webb}. In their original work they used classical optimization techniques to obtain numerical schemes, and demonstrated their performance on the linearized Euler equations. Since then, many other DRP schemes have been proposed. The main differences between existing DRP schemes is the chosen method of optimization and target/loss function. Researchers have employed methods such as simulated annealing, Lagrange multipliers, and least squares; with loss functions such as the 1-norm, 2-norm, and max-norm \cite{Zhang_DRP, Zhang_DRP2, LS_DRP, fink1, fink2}.

DRP schemes have limitations. For example, the various optimization methods often have limited wavenumebr bandwith. In this case, one usually chooses a cut-off wavenumber above which the scheme the scheme is not optimized. Other methods are designed for specific wavenumbers and while performing well on the specific ones, they struggle with approximating the solution for other wavenumbers. Hence, one of the main challenges of DRP schemes is to accurately approximate the solution in a broadband scenario.

We propose an ML based framework to create a DRP scheme, and test it in the scenario of the two-dimensional acoustic wave problem. We generate data composed of solutions of the two-dimensional wave equation. We then use a machine learning model to generate a single 
$5 \times 5$ convolutional kernel. This kernel is the FD stencil of the scheme. We achieve this by training the model using neural-network inspired methodology. Furthermore, we enrich the standard machine learning method by using a PI target function (loss). This loss function makes the algorithm ``aware" of the physical problem it is trying to solve. Using this approach we achieve improved results.	

    \section{Numerical modeling}\label{physicalProblem}

\subsection{Mathematical model}

The general formulation of the wave problem is given by:

\begin{equation} \label{waveEQ}
	\begin{cases}
		\ddot{u}(\overrightarrow{x},t)=\nabla \cdot(c^2(\overrightarrow{x})\nabla u(\overrightarrow{x},t)) & \overrightarrow{x}\in \Omega,\;t\in (0,T], \\
		u(\overrightarrow{x},0)=u_0(\overrightarrow{x}) & \overrightarrow{x}\in \Omega,\\
		\dot{u}(\overrightarrow{x},0)=v_0(\overrightarrow{x}) & \overrightarrow{x}\in \Omega,\\
		u(\overrightarrow{x},t)=f(\overrightarrow{x},t) & \overrightarrow{x}\in \partial\Omega_1,\; t\in [0,T],\\
		\frac{\partial u}{\partial \textbf{n}} (\overrightarrow{x},t)=g(\overrightarrow{x},t) & \overrightarrow{x}\in \partial\Omega_2,\; t\in [0,T], \qquad \partial\Omega_1\cup\partial\Omega_2=\partial\Omega,\\
	\end{cases}
\end{equation}
where $u(\overrightarrow{x}, t)$ is the wave amplitude or acoustic pressure, $c(\overrightarrow{x})\!>\!0$ is the wave propagation speed, $u_0$ and $v_0$ are the initial pressures and velocities respectively, and $f$ and $g$ are boundary condition of types Dirichlet and Neumann respectively. In this work we investigate the two-dimensional case. Therefore, throughout the paper $\overrightarrow{x} = (x,y)$ so that $(x, y)\in \Omega = [a_x,b_x] \times [a_y,b_y]$. We also assume that $c(\overrightarrow{x}) \equiv c$ is constant. The specific system we solve is described in section \ref{setup}. Problem \eqref{waveEQ} is well-posed and thus small changes in the problem conditions result in small changes in the solution. In addition, there exists a unique continuous solution to the problem inside the domain.

\subsection{Dispersion relation}\label{dispersion}

In this section we develop the formulae that represent the physical and numerical dispersion relations.

\subsubsection{Physical dispersion} \label{analytical_dispersion}

Given a plane-wave with wavenumber $k$, frequency $\omega$, and velocity $c$, there exists a relation $\omega \!=\! \omega(c,k)$ between them. It is called the dispersion relation, and it is dependent on the PDE and its numerical approximation. Assume that $u(x,y,t) = e^{-\textit{I}(k_xx + k_yy-\omega t)}$ satisfies the wave equation, where $\overrightarrow{k} \!=\! (k_x,k_y) \!=\! (k sin{\theta},k cos{\theta})$, and $\textit{I}$ is the imaginary unit. Inserting this into the wave equation, we get:

\begin{equation*}
	-\omega^2u = c^2(-k_x^2-k_y^2)u,
\end{equation*}
which simplifies to:
\begin{equation*}
	\omega^2 = c^2k^2.
\end{equation*}
Without loss of generality, we take into account only positive frequencies for the analysis. Hence, the expression further simplifies into the dispersion relation:	
\begin{equation}\label{dispersion_relation}
	\omega = c k.
\end{equation}	
We can express $\omega$ as a function of the wavenumber $k$:
\begin{equation*}
	\omega(k) = c(k)k.
\end{equation*}
However, since $c$ is constant we have:
\begin{equation*}
	\omega(k) = ck.
\end{equation*}	
The dispersion relation is linear in terms of $k$, which means that the PDE is non-dispersive. In physical terms, waves with different wavenumbers propagate at the same velocity.

\subsubsection{Numerical dispersion} \label{numerical_dispersion}

A general second-order in time explicit symmetric FD scheme can be formulated using convolutions as follows:

\begin{equation}\label{general_conv_drp_scheme}
	{u}_{i, j}^{n+1} = 2{u}_{i,j}^n - {u}_{i,j}^{n-1}
	+ \left(\frac{c \Delta t}{h}\right)^2 [\overrightarrow{u}^n * K]_{i,j},
\end{equation}
where $\Delta x = \Delta y = h$, $K$ is a matrix that describes the chosen stencil, and $*$ is the discrete convolution operation. We repeat the same process as in section \ref{analytical_dispersion} for \eqref{general_conv_drp_scheme} to find the numerical dispersion relation. Assume that $u(x_i,y_j,t_n) = e^{-\textit{I}(k_xi\Delta x  + k_yj\Delta y-\hat{\omega} n \Delta t)}$. Inserting into \eqref{general_conv_drp_scheme} and performing algebraic operations, we get:

\begin{equation}\label{general_kernel_dispersion}
	\alpha ^2 E*K + 4 sin^2 \left(\frac{\tilde{\omega} \Delta t}{2}\right) = 0,
\end{equation}
where $\alpha = \frac{c \Delta t}{h}$ is the CFL number, $\tilde{\omega}$ is the numerical frequency, and $E$ is a matrix of the same size as the stencil $K$, defined by $E_{i,j} = cos((i-2) h k_x + (j-2) h k_y)$.	

Extracting $\tilde{\omega}$ from \eqref{general_kernel_dispersion}, we get the ratio between the analytical frequency and the numerical frequency:

\begin{equation}\label{general_kernel_dispersion_error}
	\frac{\tilde{\omega}}{\omega} = \frac{2}{ck\Delta t} \arcsin \left(\sqrt{-\frac{\alpha}{4} ^2 E*K} \right).
\end{equation}	
The numerical dispersion relation is non-linear in terms of $k$. Thus, the FD scheme simulates a dispersive equation while the physical equation is non-dispersive.

	\section{Machine learning approach}\label{learningApproach}
We aim to find a stencil $K$ as in \eqref{general_conv_drp_scheme} of the form:

\begin{equation}\label{general_stencil}
	K = \begin{pmatrix}
		w_5 & w_4 & w_3 & w_4 & w_5  \\
		w_4 & w_2 & w_1 & w_2 & w_4  \\
		w_3 & w_1 & w_0 & w_1 & w_3  \\
		w_4 & w_2 & w_1 & w_2 & w_4  \\
		w_5 & w_4 & w_3 & w_4 & w_5  \\
	\end{pmatrix}.
\end{equation}
The coefficients $w_0, ..., w_5$ are the main result of the method we propose. The explanation for choosing this form for the stencil is explained in section \ref{constraints}.

\subsection{Data-driven problem}

We create a dataset of analytic solutions of the two-dimensional wave equation for different initial conditions. Specifically, we fix the grid size, velocity, and maximal wavenumber of the initial conditions while only changing the wave amplitudes. We choose a wavenumber for which classical numerical methods exhibit high numerical dispersion as the maximal wavenumber in the dataset. The dataset includes solutions produced from lower wavenumbers as well.

The dataset is then divided into samples and labels. The samples corresponding to each initial condition in the dataset are given by: $\{(u^{n-1}, u^n)\}_{n=1}^{N_t - 1}$. Similarly, the corresponding labels are given by $\{(u^{n+1})\}_{n=1}^{N_t - 1}$. A detailed construction of the dataset is given in section \ref{setup}.

The input to the network is a set of samples $u^{n-1}$ and $u^n$ where both are $N_x X N_y$ matrices.
We then apply the convolutional kernel exclusively on $u^n$ and predict $u^{n+1}$ via the procedure shown in \eqref{general_conv_drp_scheme}
We use a mean squared error (MSE) loss function to proceed with the kernel optimization. It is defined by:  \\
\begin{equation}\label{MSE}
	MSE(\hat{u}^{n+1}, u^{n+1}) = \frac{1}{N_x N_y}\sum\limits_{i=1}^{N_x}\sum\limits_{j=1}^{N_y}\left(\hat{u}^{n+1}_{i,j} - {u}^{n+1}_{i,j}\right)^2,
\end{equation}
where $\hat{u}^{n+1}$ is the network's prediction and ${u}^{n+1}$ is the exact solution.

\subsection{PI iterative loss} \label{phyloss}

We seek to improve the capabilities of the model by adding a PI loss to the training process. In this case, we use an iterative loss. The iterative loss computes three time marching steps using the predictions of the model, during training. These three steps produce the wave pressure. We compare the pressures to the ones taken from the analytical solution of the problem at the same time steps. If the predictions are accurate, the pressures are similar and the loss is low. This loss acts as a penalty term for the model and, during training, enhances the convergence of the model. The iterative loss is composed of three different loss terms:

\underline{Regular MSE}:
$MSE(\hat{u}^{(n+1)}_{pred}, u^{(n+1)}_{true})$, where:
\begin{equation*}
	\hat{u}^{(n+1)}_{pred} = \alpha ^2 (u^{(n)}_{true} * K) + 2u^{(n)}_{true} - u^{(n-1)}_{true},
\end{equation*}

\underline{Semi-inferred MSE}:
$MSE(\hat{u}^{(n+2)}_{pred}, u^{(n+2)}_{true})$, where:
\begin{equation*}
	\hat{u}^{(n+2)}_{pred} = \alpha ^2 (u^{(n+1)}_{pred} * K) + 2u^{(n+1)}_{pred} - u^{(n)}_{true},
\end{equation*}

\underline{Fully-inferred MSE}:
$MSE(\hat{u}^{(n+3)}_{pred}, u^{(n+3)}_{true})$, where:
\begin{equation*}
	\hat{u}^{(n+3)}_{pred} = \alpha ^2 (u^{(n+2)}_{pred} * K) + 2u^{(n+2)}_{pred} - u^{(n+1)}_{pred},
\end{equation*}
Finally, the iterative loss is defined by:
\begin{equation}\label{iterative_loss_definition}
	L_{iterative} = MSE(\hat{u}^{(n+1)}_{pred}, u^{(n+1)}_{true}) + MSE(\hat{u}^{(n+2)}_{pred}, u^{(n+2)}_{true}) + MSE(\hat{u}^{(n+3)}_{pred}, u^{(n+3)}_{true}),
\end{equation}

\subsection{PI constraints} \label{constraints}
In addition to the PI loss, we added constraints on the stencil $K$ \eqref{general_stencil} during the training process:

\begin{itemize}
	\item Consistency: we require the sum of all elements to be zero to enforce consistency. Using Taylor expansion, this ensures at least first order accuracy. Mathematically:
	\begin{equation*}
		w_0 + 4w_1 + 4w_2 + 4w_3 + 8w_4 + 4w_5 = 0
	\end{equation*}
	\item Symmetry: the stencil must be symmetric with regards to its central point. The symmetry is shown in the definition of the stencil \eqref{general_stencil}. This cancels out terms in the Taylor expansion, and combined with the following equation ensures at least second order accuracy \cite{high_order}:
	\begin{equation*}
		w_1 + 2w_2 + 4w_3 + 10w_4 + 8w_5 = 1
	\end{equation*}	
	\item Fourth order: we add a constraint derived from the Taylor theorem that guarantees fourth order accuracy \cite{high_order}:
	\begin{equation*}
		\begin{cases}
			w_2 + 8w_4 + 16w_5 = 0 \\
			\frac{1}{12}w_1 + \frac{1}{6}w_2 + \frac{4}{3}w_3 + \frac{17}{6}w_4 + \frac{8}{3}w_5 = 0
		\end{cases}
	\end{equation*}		
\end{itemize}

These three constraints together define the rules for our fourth order stencil $K_4$. We also experimented with a second order stencil $K_2$ by only applying the consistency and symmetry constraints, and discarding unnecessary weights. The second order stencil is defined as:

 	\begin{equation}\label{k2}
 	K_2 = \begin{pmatrix}
 		0 & 0 & w_3 & 0 & 0  \\
 		0 & 0 & w_1 & 0 & 0  \\
 		w_3 & w_1 & w_0 & w_1 & w_3  \\
 		0 & 0 & w_1 & 0 & 0  \\
 		0 & 0 & w_3 & 0 & 0  \\
 	\end{pmatrix}.
 \end{equation}

	\section{Numerical tests and results} \label{tests_and_results}

\subsection{Experiment setup} \label{setup}

The specific equation being solved throughout this section is given by:
\begin{equation} \label{specificWaveEQ2d}
	\begin{cases}
		\ddot{u}(x, y, t)=c^2 (u_{xx}(x, y, t) + u_{yy}(x, y, t)) & 0\leq x,y\leq L,\; 0\leq t \leq T, \\
		u(x, y, 0)=u_0(x, y) & 0\leq x,y\leq L, \\
		\dot{u}(x, y, 0)=v_0(x,y) & 0\leq x,y\leq L, \\
		u(0, y, t)=u(L, y, t)=0 & 0\leq y\leq L,\; 0\leq t \leq T, \\
		u(x, 0, t)=u(x, L, t)=0 & 0\leq x\leq L,\; 0\leq t \leq T. \\
	\end{cases}
\end{equation}

The model was trained on solutions of this equation corresponding to sine initial condition. A generic sample from this dataset is the projection of the following solution onto a numerical grid:
\begin{equation}\label{general_solution_2d}
	\sum_{n=1}^{d}{a_n sin\left(\frac{\pi n}{L} x\right)sin\left(\frac{\pi n}{L} y\right) cos\left(\sqrt{2}\frac{\pi n}{L} c t\right)},
\end{equation}
where d is the degree of the specified trigonometric polynomial.

\begin{figure}[!h]
	\centering
	\includegraphics[width=16cm]{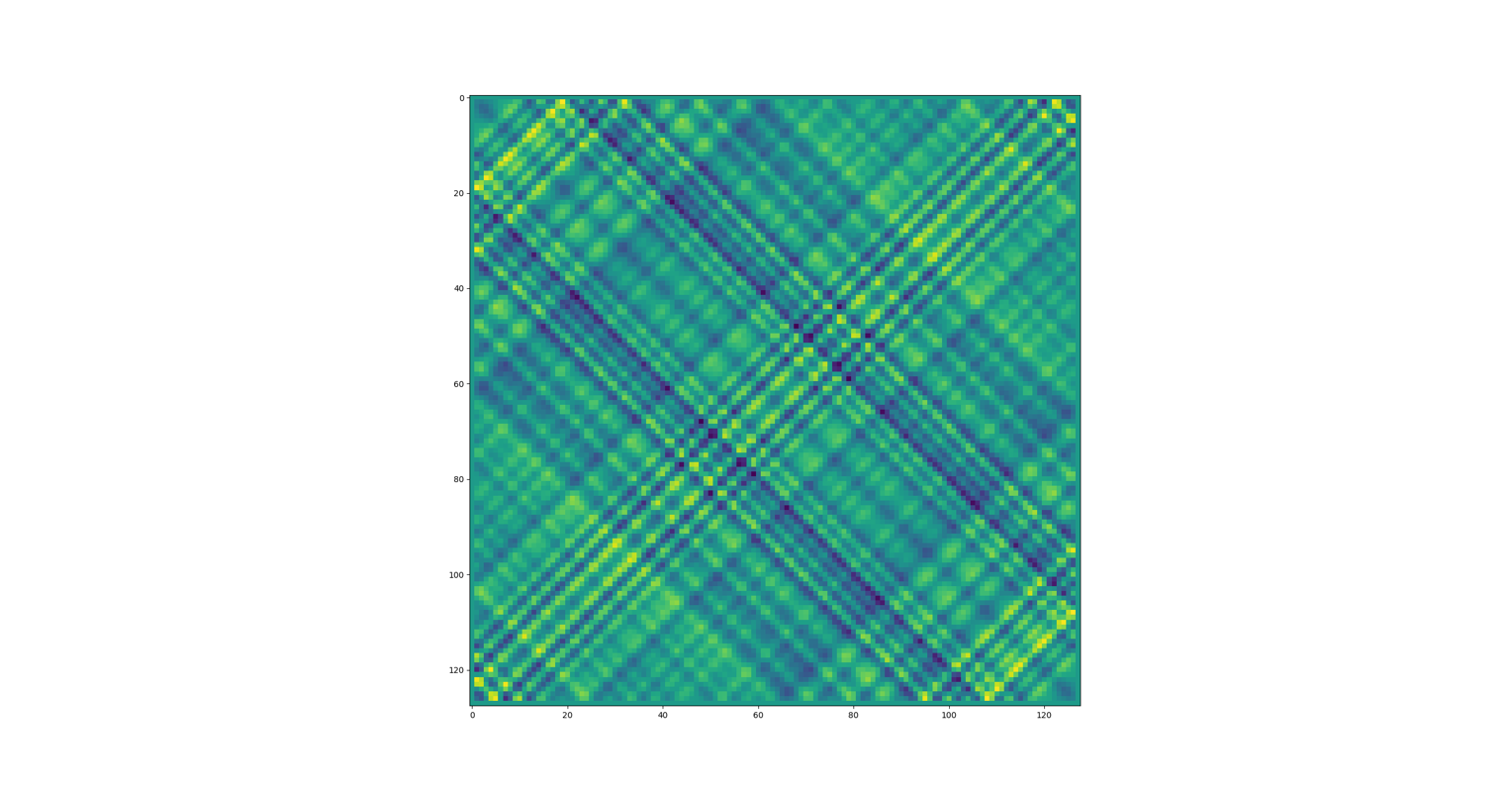}
	\caption{A trigonometric polynomial of degree 80 $\pi$.}
	\label{fig:example_sine}
\end{figure}

The following physical parameters were used during the data generation process: $L \!=\! 1$, $T \!=\! 0.1$, $c\!=5\!$, and $v_0(x)\!=0\!$. We use solutions of the form \eqref{general_solution_2d} where the maximal degree of the trigonometric polynomial was $80 \pi$. One such random trigonometric polynomial is shown in figure \ref{fig:example_sine}. These solutions were discretized on a spatial grid of size $128 \times 128$ using $400$ time steps. We generated 12 such solutions which we then split into training, testing, and validation sets according to standard ML practice.

We compared the results of the following schemes: classical second order central differences (Classic2), classical fourth order central differences (Classic4), and an existing optimized DRP scheme (DRP) \cite{Zhang_DRP}.

\

\subsection{Results}
After training the model with the physically-informed components, the following second and fourth order stencils were obtained:

\begin{equation}\label{drp_kernel4}
	K_2 = \begin{pmatrix}
         0 & 0 & -0.10458 & 0 &  0 \\
         0 & 0 &  1.41831 & 0 &  0 \\
        -0.10458 & 1.41831 & -5.25494 & 1.41831 & -0.10458 \\
         0 & 0 &  1.41831 & 0 &  0 \\
         0 & 0 & -0.10458 & 0 &  0
		
	\end{pmatrix},	
\end{equation}
	
\begin{equation}\label{drp_kernel4}
	K_4 = \begin{pmatrix}
        -0.00917 &  0.04473 & -0.15445 &  0.04473 & -0.00917 \\
         0.04473 & -0.21106 &  1.66600 & -0.21106 &  0.04473 \\
        -0.15445 &  1.66600 & -5.52311 &  1.66600 & -0.15445 \\
         0.04473 & -0.21106 &  1.66600 & -0.21106 &  0.04473 \\
        -0.00917 &  0.04473 & -0.15445 &  0.04473 & -0.00917
		
	\end{pmatrix},	
\end{equation}
We tested it over several scenarios by changing the initial condition. For the purpose of error analysis we use the $L^2$ grid norm defined as:

\begin{equation}
	|| \cdot ||_{2, h} = \sqrt{h^2 \Delta t \sum\limits_{n=1}^{N_t}\sum\limits_{i=1}^{N_x}\sum\limits_{j=1}^{N_y}\left|u_{exact}(x_i, y_j, t_n) - {u}^{n}_{i,j}\right|^2},
\end{equation}

assuming $\Delta x = \Delta y = h$.

Similarly, the convergence rate for two discretizations $h_1$ and $h_2$ is given by:

\begin{equation}
	Rate =
	\frac{{\log \left( \frac{|| \cdot ||_{2, h_1}}{|| \cdot ||_{2, h_2}} \right)}}{{\log \left( \frac{h_1}{h_2} \right) }}
\end{equation}

Note that like other fourth order methods, these stencils require a special treatment near the boundary. We used Fornberg's \cite{Fornberg} method to generate one-sided FD approximations at the points near the boundary. To keep the order of accuracy consistent these approximations were of the same order as the used stencil in each case.

\begin{figure}[!h]
	\centering
	\includegraphics[width=16cm]{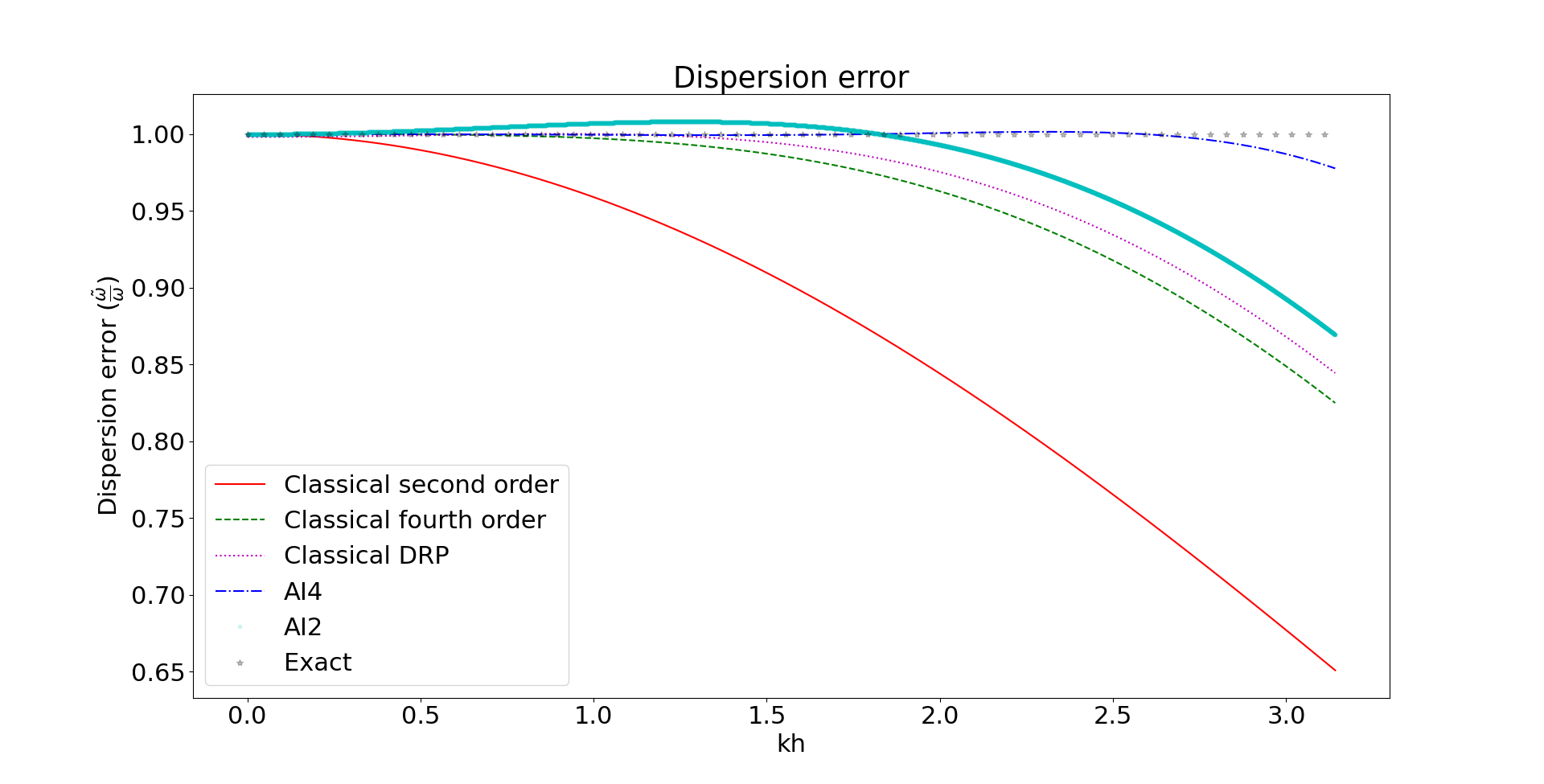}
	\caption{Dispersion error over $kh$ - the wavenumber multiplied by the spatial discretization parameter, where 1 is a perfect match between the physical and numerical frequencies.}
	\label{fig:dispersion_full_range}
\end{figure}

We plot the numerical dispersion error of the various methods in figure \ref{fig:dispersion_full_range}. It is clear to see that the standard second and fourth order schemes achieve the lowest accuracy range, while the various DRP schemes outperform them. More specifically, the fourth order AI scheme stays the most accurate for higher wavenumbers.

\

\subsubsection{Experiment I - low wavenumber} \label{scenario1}
We begin by testing the scheme on a single sine initial condition with a low wavenumber:

\begin{equation*}
	u(x, y, 0) = sin(\pi x) sin(\pi y)
\end{equation*}
with $c=1, T=0.08, L=1$.

The analytical solution is:

\begin{equation*}
	u(x, y, t) = sin(\pi x) sin(\pi y) cos(\sqrt{2}\pi t)
\end{equation*}

As seen in table \eqref{table:errors_scenario1}, the proposed fourth order scheme (AI4) outperforms all other methods for all discretization parameters. The proposed second order scheme (AI2) yields a lower error than the classical second order FD (Classic2) and the benchmark DRP scheme. In table \eqref{table:convergence_scenario1} we can see that the computed convergence rate of the proposed stencils agrees with the theoretical one.

\begin{table}[h!]
	\begin{center}\small
		\begin{tabular}{||c c | c c c c c||}
			\hline
			$N_x \times N_y$ & $N_t$ & Classic2 & Classic4 & DRP & AI4 & AI2 \\ [0.5ex]
			\hline\hline
			$10 \times 10$ & 25 & 3.964201e-05 & 1.341580e-06 & 3.108537e-06 & \textbf{1.247523e-06} & 8.963064e-06 \\
			$20 \times 20$ & 100 & 8.950977e-06 & 3.337217e-08 & 6.122411e-06 & \textbf{1.973002e-08} & 2.227419e-06 \\
			$40 \times 40$ & 400 & 2.127100e-06 & 1.580853e-09 & 6.961017e-06 & \textbf{3.386081e-10} &          5.394590e-07 \\
			$80 \times 80$ & 1600 & 5.185358e-07 & 9.298174e-11 & 7.162837e-06 & \textbf{1.402962e-11} &          1.320516e-07 \\
			$160 \times 160$ & 6400 &           1.280114e-07 &           5.648513e-12 &      7.211200e-06 &    \textbf{8.376164e-13} &          3.262916e-08 \\
			\hline
		\end{tabular}
		\caption{Comparison of $L^2$ errors for experiment \ref{scenario1}.}
		\label{table:errors_scenario1}
	\end{center}
\end{table}

\begin{table}[h!]
	\begin{center}\small
		\begin{tabular}{||c c | c c c c c||}
			\hline
			$N_x \times N_y$ & $N_t$ & Classic2 & Classic4 & DRP & AI4 & AI2 \\ [0.5ex]
			\hline\hline
			$20 \times 20$ & 100 & 2.14691 & 5.32914 & N/A & 5.98253 & 2.00862 \\
			$40 \times 40$ & 400 & 2.07316 & 4.39987 & N/A & 5.86463 & 2.04579 \\
			$80 \times 80$ & 1600 & 2.03637 & 4.08761 &  N/A & 4.59307 & 2.03041 \\
			$160 \times 160$ & 6400 & 2.01817 & 4.04100 & N/A & 4.06604 & 2.01687 \\
			\hline
		\end{tabular}
		\caption{Comparison of $L^2$ convergence rates for experiment \ref{scenario1}.}
		\label{table:convergence_scenario1}
	\end{center}
\end{table}

\begin{figure}[!h]
	\centering
	\includegraphics[width=16cm]{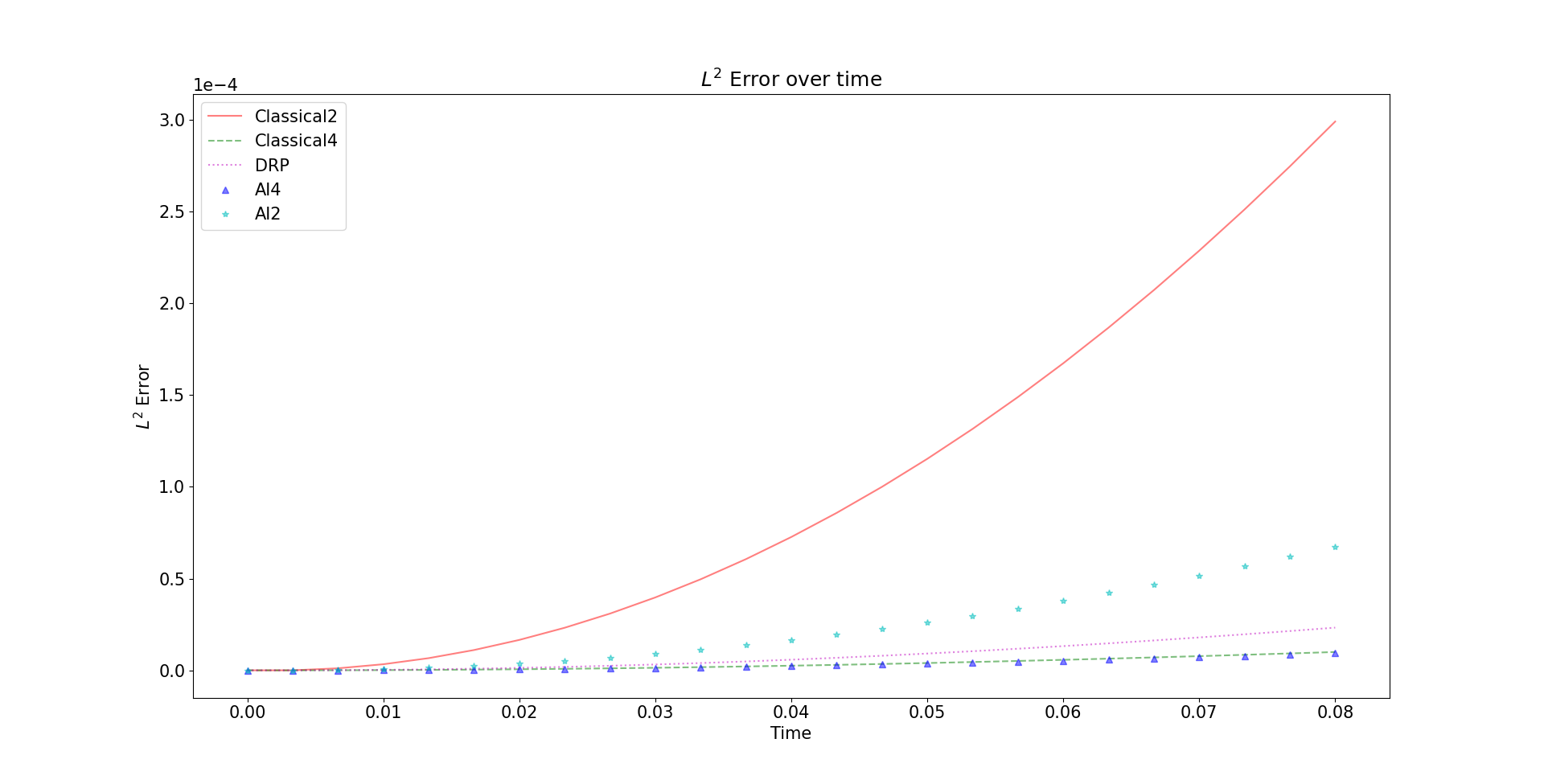}
	\caption{Error over time comparing the various numerical methods for experiment \ref{scenario1} for the coarsest grid.}
	\label{fig:error_over_time_1}
\end{figure}\

\subsubsection{Experiment II - higher wavenumber} \label{scenario2}
We use the same parameters as in \eqref{scenario1}, but we increase the wavenumber of the initial condition:

\begin{equation*}
	u(x, y, 0) = sin(20\pi x) sin(20\pi y)
\end{equation*}
with $c=1, T=0.08, L=1$.

The analytical solution is:

\begin{equation*}
	u(x, y, t) = sin(20\pi x) sin(20\pi y) cos(\sqrt{40}\pi t)
\end{equation*}

In this case (\eqref{table:errors_scenario2}) we see once again that AI4 achieves the lowest errors amongst the compared methods. We see an especially significant drop in the last final discretization where AI4 outperforms all other methods by a factor of one to three orders of magnitudes. Interestingly, the AI2 is slightly more accurate than the Classic4 and benchmark DRP for very coarse grids. However, when we refine the grid the latter outperform the former. Consequently, in this scenario the AI2 stencil is sufficient for extremely coarse grids, but the AI4 should be used for more refined grids. This is backed up by the convergence rates given in \eqref{table:convergence_scenario2}.

\begin{table}[h!]
\begin{center}\small
	\begin{tabular}{||c c | c c c c c||}
		\hline
		$N_x \times N_y$ & $N_t$ & Classic2 & Classic4 & DRP & AI4 & AI2 \\ [0.5ex]
		\hline\hline
		$10 \times 10$ & 25 & 8.119670e-02 & 7.858965e-02 & 7.809909e-02 & \textbf{7.470846e-02} & 7.751773e-02 \\
		$20 \times 20$ & 100 & 2.093424e-02 & 1.665968e-02 & 1.623471e-02 & \textbf{1.483528e-02} & 1.576807e-02 \\
		$40 \times 40$ & 400 & 2.024183e-02 & 5.065186e-03 & 3.952975e-03 & \textbf{2.638155e-03} & 2.865844e-03 \\
		$80 \times 80$ & 1600 & 4.941423e-03 & 4.480877e-04 & 3.282980e-04 & \textbf{3.228985e-04} & 8.576374e-04 \\
		$160 \times 160$ & 6400 & 1.222075e-03 & 2.154743e-05 & 1.222193e-04 & \textbf{6.101200e-06} & 2.810811e-04  \\
		\hline
	\end{tabular}
	\caption{Comparison of $L^2$ errors for experiment \ref{scenario2}.}
	\label{table:errors_scenario2}
\end{center}
\end{table}

\begin{table}[h!]
	\begin{center}\small
		\begin{tabular}{||c c | c c c c c||}
			\hline
			$N_x \times Ny$ & $Nt$ & Classic2 & Classic4 & DRP & AI4 & AI2 \\ [0.5ex]
			\hline\hline
			$20 \times 20$ & 100 & 1.95556 & 2.23798 &           2.26622 &         2.33224 &               2.29752 \\
			$40 \times 40$ & 400 & 0.04852 &                1.71767 &           2.03807 &         2.49143 &               2.45997 \\
			$80 \times 80$ & 1600 & 2.03434 &                3.49876 &           3.58986 &         3.03038 &               1.74052 \\
			$160 \times 160$ & 6400 & 2.01559 &                4.37819 &           1.42553 &         5.72584 &               1.60938 \\
			\hline
		\end{tabular}
		\caption{Comparison of $L^2$ convergence rates for experiment \ref{scenario2}.}
		\label{table:convergence_scenario2}
	\end{center}
\end{table}

\begin{figure}[!h]
	\centering
	\includegraphics[width=16cm]{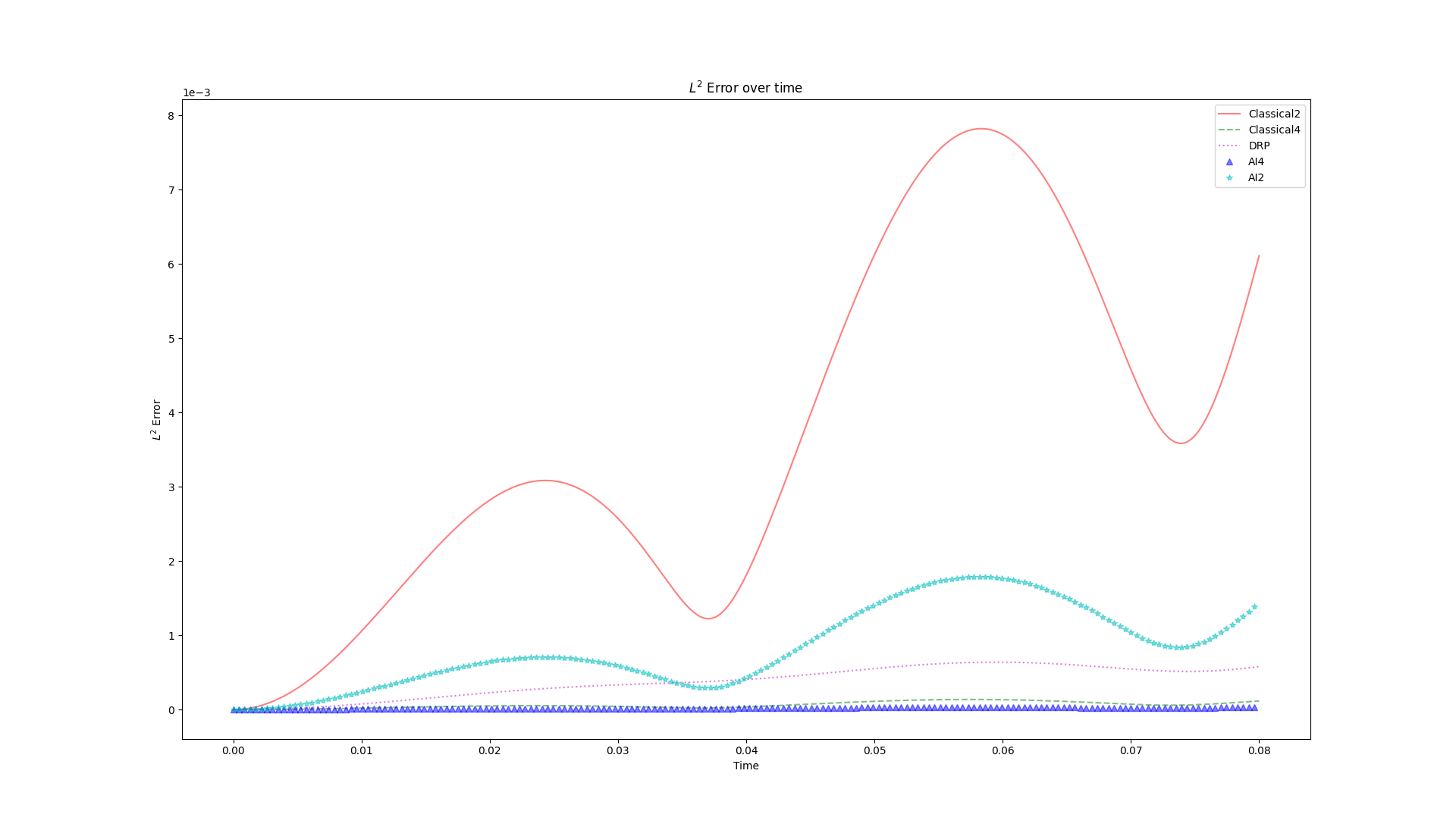}
	\caption{Error over time comparing the various numerical methods for experiment \ref{scenario2} for the most refined grid.}
	\label{fig:error_over_time_20}
\end{figure}\

\subsubsection{Experiment III} \label{scenario3}
Using sines as basis functions for training and testing the model raises the question whether the model can
perform well given an initial condition that was not formed using a finite linear combination of sine functions.
Hence, we test it on a Ricker wavelet \cite{ricker1}, \cite{ricker2}. The Ricker wavelet is very common in numerical modeling of seismic applications. It is defined by:
\begin{equation}\label{ricker_definition}
	f(t) = (1 - 2\pi^2 \sigma ^2 t^2)e^{-\pi^2 \sigma ^2 t^2}	
\end{equation}
for a given frequency $\sigma$. We use the Ricker wavelet as a point source term for the wave equation \eqref{specificWaveEQ2d}. In order to simulate a realistic use-case, we set $c = 1,500, T=0.5, L=2,000$, corresponding to a homogeneous domain of size $2km \times 2km$. We locate a Ricker source with frequency $\sigma = 20$Hz at the center of the domain. As for discretization parameters we choose a coarse grid for the given frequency by setting $N_t = 500$ and $N_x = 240$.

A closed-form analytic solution is not easily obtained for this problem. Hence, all solutions were compared to a highly accurate 36th order finite difference solution devised by Devito \cite{devito1}, \cite{devito2}.

\begin{figure}[!h]
	\centering
	\includegraphics[width=12cm]{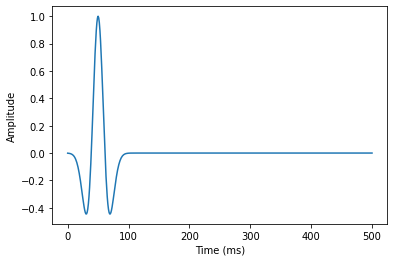}
	\caption{One-dimensional cut of Ricker wavelet.}
	\label{fig:ricker_ic_1d}
\end{figure}

\begin{figure}[!h]
	\centering
	\includegraphics[scale=0.38]{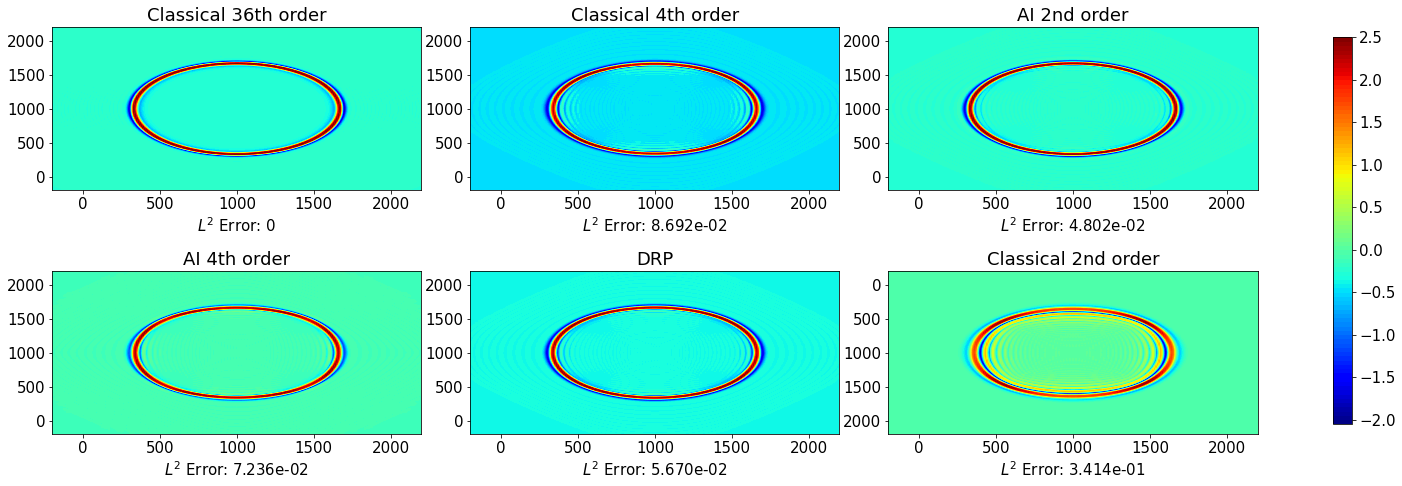}
	\caption{Comparison of six different numerical method corresponding to scenario \eqref{scenario3} at t=0.5.}
	\label{fig:ricker_solutions}
\end{figure}

\begin{figure}[!h]
	\centering
	\includegraphics[scale=0.31]{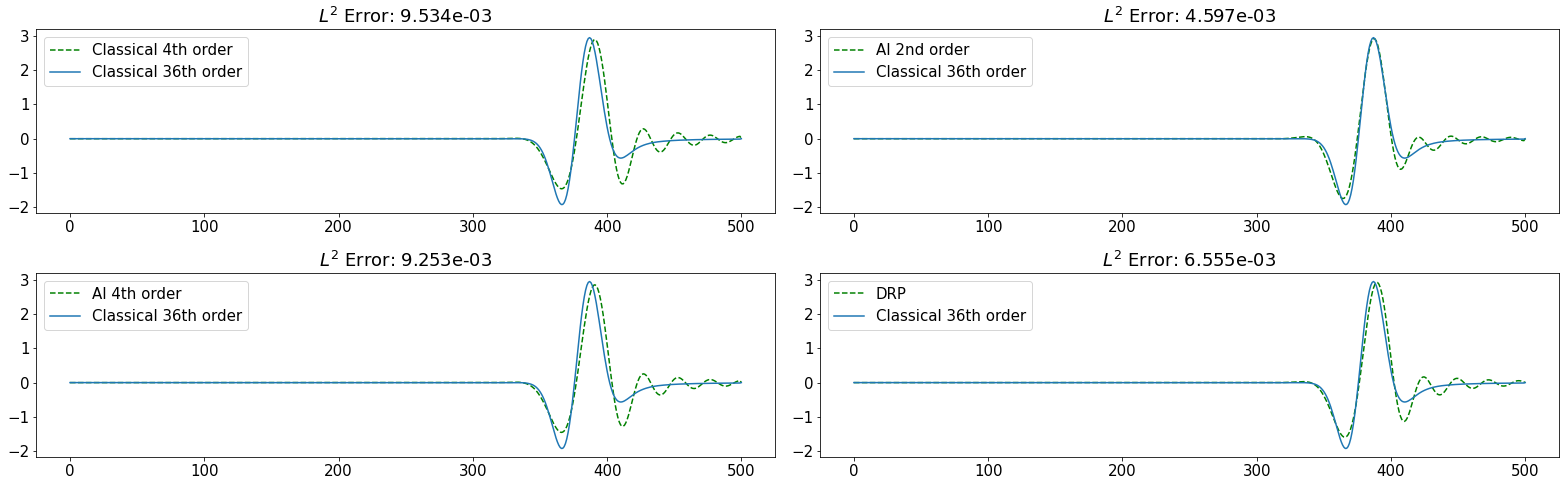}
	\caption{Comparison of waveforms corresponding to scenario \eqref{scenario3} obtained through a reciever placed at $(200m, 500m)$.}
	\label{fig:ricker_recievers}
\end{figure}

\newpage

A comparison of classical Taylor-based FD methods, DRP, and our proposed method for this experiment can be seen in figure \eqref{fig:ricker_solutions}. The second order AI scheme achieves the best accuracy in comparison to the reference solution. The classical second and fourth order schemes exhibit visible numerical dispersion. A similar comparison under the same setting is shown in \eqref{fig:ricker_recievers}, where waveforms obtained by recievers are shown. Hence, our scheme performs well for Ricker wavelets as well and not just for trigonometric polynomials.

	\section{Conclusion}\label{Conclusion}
In this work we introduced an FD method to accurately solve PDEs in the setting of coarse grids and high wavenumbers. We achieved this by developing a DRP-like scheme with optimized FD coefficients. We found the stencil by using an ML convolutional approach. We added a physically-informed component which helped us obtain accurate results. Finally, we tested our model on different initial conditions with a variety of discretizations, and saw that the proposed method performed well and converged according to theory.

We plan on expanding the model to other difficult problems in different domains. For example, Maxwell's equations pose an interesting challenge due to the high propagation velocity involved in electromagnetism. A natural expansion would be to conduct a similar process in the case of non-constant wave velocity. Another future research direction would be to solve similar problems with high wavenumbers in 3D, and find a way to generalize the proposed method.  
	  
    \bibliographystyle{unsrtnat}
    \bibliography{references}
	
\end{document}